%
%
%
\input amssym.tex
\newif\ifsect\newif\iffinal
\secttrue\finalfalse
\def\thm #1: #2{\medbreak\noindent{\bf #1:}\if(#2\thmp\else\thmn#2\fi}
\def\thmp #1) { (#1)\thmn{}}
\def\thmn#1#2\par{\enspace{\sl #1#2}\par
        \ifdim\lastskip<\medskipamount \removelastskip\penalty 55\medskip\fi}
\def\square{{\msam\char"03}}
\def\qedn{\thinspace\null\nobreak\hfill\square\par\medbreak}
\def\pf{\ifdim\lastskip<\smallskipamount \removelastskip\smallskip\fi
        \noindent{\sl Proof\/}:\enspace}

%
%
%
%
\newcount\parano
\newcount\eqnumbo
\newcount\thmno
\newcount\versiono
\newcount\remno
\newbox\notaautore
\def\neweqt#1$${\xdef #1{(\number\parano.\number\eqnumbo)}
    \eqno #1$$
    \global \advance \eqnumbo by 1}
\def\newrem #1\par{\global \advance \remno by 1
    \medbreak
{\bf Remark \the\remno:}\enspace #1\par
\ifdim\lastskip<\medskipamount \removelastskip\penalty 55\medskip\fi}
\def\newthmt#1 #2: #3{\xdef #2{\number\parano.\number\thmno}
    \global \advance \thmno by 1
    \medbreak\noindent
    {\bf #1 #2:}\if(#3\thmp\else\thmn#3\fi}
\def\neweqf#1$${\xdef #1{(\number\eqnumbo)}
    \eqno #1$$
    \global \advance \eqnumbo by 1}
\def\newthmf#1 #2: #3{\xdef #2{\number\thmno}
    \global \advance \thmno by 1
    \medbreak\noindent
    {\bf #1 #2:}\if(#3\thmp\else\thmn#3\fi}
\def\forclose#1{\hfil\llap{$#1$}\hfilneg}
\def\newforclose#1{
	\ifsect\xdef #1{(\number\parano.\number\eqnumbo)}\else
	\xdef #1{(\number\eqnumbo)}\fi
	\hfil\llap{$#1$}\hfilneg
	\global \advance \eqnumbo by 1
	\iffinal\else\rsimb#1\fi}
\def\forevery#1#2$${\displaylines{\let\eqno=\forclose
        \let\neweq=\newforclose\hfilneg\rlap{$\qquad\quad\forall#1$}\hfil#2\cr}$$}
\def\noNota #1\par{}
\def\today{\ifcase\month\or
   January\or February\or March\or April\or May\or June\or July\or August\or
   September\or October\or November\or December\fi
   \space\number\year}
\def\inizia{\ifsect\let\neweq=\neweqt\else\let\neweq=\neweqf\fi
\ifsect\let\newthm=\newthmt\else\let\newthm=\newthmf\fi}
\def\bititolo{\empty}
\gdef\begin #1 #2\par{\xdef\titolo{#2}
\ifsect\let\neweq=\neweqt\else\let\neweq=\neweqf\fi
\ifsect\let\newthm=\newthmt\else\let\newthm=\newthmf\fi
\iffinal\let\Nota=\noNota\fi
\centerline{\titlefont\titolo}
\if\bititolo\empty\else\medskip\centerline{\titlefont\bititolo}
\xdef\titolo{\titolo\ \bititolo}\fi
\bigskip
\centerline{\bigfont
\autore \ifvoid\notaautore\else\footnote{${}^1$}{\unhbox\notaautore}\fi}
\bigskip\if\istituto!\centerline{\today}\else
\centerline{\istituto}
\centerline{\indirizzo}
\centerline{\email}
\medskip
\centerline{#1\ \anno}\fi
\bigskip\bigskip
\ifsect\else\global\thmno=1\global\eqnumbo=1\fi}
\def\anno{2011}
\def\raggedleft{\leftskip2cm plus1fill \spaceskip.3333em \xspaceskip.5em
\parindent=0pt\relax}
\font\titlefont=cmssbx10 scaled \magstep1
\font\bigfont=cmr12
\font\eightrm=cmr8
\font\sc=cmcsc10
\font\bbr=msbm10
\font\sbbr=msbm7
\font\ssbbr=msbm5
\font\msam=msam10

\font\bfm=cmmib10

\def\ca #1{{\cal #1}}

\nopagenumbers
\binoppenalty=10000
\relpenalty=10000
\newfam\amsfam
\textfont\amsfam=\bbr \scriptfont\amsfam=\sbbr \scriptscriptfont\amsfam=\ssbbr
\newfam\boldifam
\textfont\boldifam=\bfm
\let\de=\partial
\let\eps=\varepsilon
\let\phe=\varphi
\def\Hol{\mathop{\rm Hol}\nolimits}

\def\Fix{\mathop{\rm Fix}\nolimits}
\def\Sing{\mathop{\rm Sing}\nolimits}
\def\End{\mathop{\rm End}\nolimits}

\def\Res{\mathop{\rm Res}\nolimits}
\def\Re{\mathop{\rm Re}\nolimits}

\def\Jac{\mathop{\rm Jac}\nolimits}

\def\bigoperp{\mathop{\hbox{$\bigcirc\kern-11.8pt\perp$}}\limits}

\mathchardef\void="083F
\mathchardef\ellb="0960
\mathchardef\taub="091C
\def\C{{\mathchoice{\hbox{\bbr C}}{\hbox{\bbr C}}{\hbox{\sbbr C}}
{\hbox{\sbbr C}}}}

\def\N{{\mathchoice{\hbox{\bbr N}}{\hbox{\bbr N}}{\hbox{\sbbr N}}
{\hbox{\sbbr N}}}}
\def\P{{\mathchoice{\hbox{\bbr P}}{\hbox{\bbr P}}{\hbox{\sbbr P}}
{\hbox{\sbbr P}}}}

\newcount\notitle
\notitle=1
\headline={\ifodd\pageno\rhead\else\lhead\fi}
\def\rhead{\ifnum\pageno=\notitle\iffinal\hfill\else\hfill\tt Version
\the\versiono; \the\day/\the\month/\the\year\fi\else\hfill\eightrm\titolo\hfill
\folio\fi}
\def\lhead{\ifnum\pageno=\notitle\hfill\else\eightrm\folio\hfill\autore\hfill
\fi}
\newbox\bibliobox
\def\setref #1{\setbox\bibliobox=\hbox{[#1]\enspace}
    \parindent=\wd\bibliobox}
\def\biblap#1{\noindent\hang\rlap{[#1]\enspace}\indent\ignorespaces}
\def\art#1 #2: #3! #4! #5 #6 #7-#8 \par{\biblap{#1}#2: {\sl #3\/}.
    #4 {\bf #5} (#6)\if.#7\else, \hbox{#7--#8}\fi.\par\smallskip}
\def\book#1 #2: #3! #4 \par{\biblap{#1}#2: {\bf #3.} #4.\par\smallskip}
\def\coll#1 #2: #3! #4! #5 \par{\biblap{#1}#2: {\sl #3\/}. In {\bf #4,}
#5.\par\smallskip}
\def\pre#1 #2: #3! #4! #5 \par{\biblap{#1}#2: {\sl #3\/}. #4, #5.\par\smallskip}
\def\Bittersweet{\pdfsetcolor{0. 0.75 1. 0.24}} 
\def\Black{\pdfsetcolor{0. 0. 0. 1.}} 
\def\raggedleft{\leftskip2cm plus1fill \spaceskip.3333em \xspaceskip.5em
\parindent=0pt\relax}
\def\Nota #1\par{\medbreak\begingroup\Bittersweet\raggedleft
#1\par\endgroup\Black
\ifdim\lastskip<\medskipamount \removelastskip\penalty 55\medskip\fi}
\newcount\defno\newcount\exno
\def\smallsect #1. #2\par{\bigbreak\noindent{\bf #1.}\enspace{\bf #2}\par
    \global\parano=#1\global\eqnumbo=1\global\thmno=1\global\defno=0\global\remno=0
    \nobreak\smallskip\nobreak\noindent\message{#2}}
\def\newdef #1\par{\global \advance \defno by 1
    \medbreak
{\bf Definition \the\defno:}\enspace #1\par
\ifdim\lastskip<\medskipamount \removelastskip\penalty 55\medskip\fi}
\def\newex #1\par{\global \advance \exno by 1
    \medbreak
{\bf Example \the\exno:}\enspace #1\par
\ifdim\lastskip<\medskipamount \removelastskip\penalty 55\medskip\fi}
\finaltrue
\versiono=1
%



\sectfalse

\def\autore{Marco Abate\footnote{}{\eightrm 2010 Mathematics Subject Classification: Primary 32H04; 32H50, 14M20, 37F75, 37F99, 58J20}}
\def\indirizzo{Largo Pontecorvo 5, 56127 Pisa, Italy.}
\def\istituto{Dipartimento di Matematica, Universit\`a di Pisa}
\def\email{E-mail: abate@dm.unipi.it}

\begin {June} Index theorems for meromorphic self-maps of the projective space

\centerline{\sl Dedicated to Jack Milnor for his 80th birthday}
\bigskip
{{\narrower{\sc Abstract.} In this short note we would like to show how it is possible to use
techniques introduced in the theory of local dynamics of holomorphic germs tangent to the
identity to study global meromorphic self-maps of the complex projective space. In particular
we shall show how a meromorphic map $f\colon\P^n\dasharrow\P^n$ induces a
holomorphic foliation of $\P^n$ in Riemann surfaces, whose singular points are exactly
the fixed points and the indeterminacy points of~$f$; and we shall prove three index 
theorems, relating suitably defined local residues at the fixed and indeterminacy points
of $f$ with Chern classes of $\P^n$.

}}

\bigskip\noindent
In this short note we would like to show how the techniques introduced in [ABT1] (see also
[ABT2, 3], and [AT2]) for studying the local dynamics of holomorphic germs tangent to the identity can be used to study global meromorphic self-maps
of the complex projective space~$\P^n$. More precisely, we shall prove the following index theorem:

\newthm Theorem \index: Let $f\colon\P^n\dasharrow\P^n$ be a meromorphic self-map of degree~$\nu+1\ge 2$ of the complex $n$-dimensional projective space. Let $\Sigma(f)=\Fix(f)\cup I(f)$ be
the union of the indeterminacy set $I(f)$ of~$f$ and the fixed points set $\Fix(f)$ of~$f$.
Let $\Sigma(f)=\sqcup_\alpha \Sigma_\alpha$ be the decompositon of~$\Sigma$ in
connected components, and denote by $N$ the tautological line bundle of $\P^n$. Then:
\smallskip
\item{\rm(i)} we can associate to each $\Sigma_\alpha$ a complex number $\Res^1(f,\Sigma_\alpha)\in\C$, depending only on the local behavior of $f$ nearby~$\Sigma_\alpha$, so that
$$
\sum_\alpha \Res^1(f,\Sigma_\alpha)=\int_{\P^n}c_1(N)^n=(-1)^n\;;
$$
\item{\rm(ii)} given a homogeneous symmetric polynomial $\phe\in\C[z_1,\ldots,z_n]$ of degree~$n$
we can associate to each $\Sigma_\alpha$ a complex number $\Res^2_\phe(f,\Sigma_\alpha)\in\C$, depending only on the local behavior of $f$ nearby~$\Sigma_\alpha$,
such that
$$
\sum_\alpha \Res^2_\phe(f,\Sigma_\alpha)=\int_{\P^n}\phe(T\P^n- N^{\otimes\nu})\;;
$$
\item{\rm(iii)} if $\nu>1$, given a homogeneous symmetric polynomial $\psi\in\C[z_0,\ldots,z_n]$ of degree~$n$
we can associate to each $\Sigma_\alpha$ a complex number $\Res^3_\psi(f,\Sigma_\alpha)\in\C$, depending only on the local behavior of $f$ nearby~$\Sigma_\alpha$,
such that
$$
\sum_\alpha \Res^3_\psi(f,\Sigma_\alpha)=\int_{P^n}\psi\bigl((T\P^n\oplus N)- N^{\otimes\nu}\bigr)\;.
$$

In this statement, if $\phe\in\C[z_1,\ldots,z_r]$ is a homogeneous symmetric polynomial and $E$
and $F$ are vector bundles over $\P^n$, we put
$$
\phe(E-F)=\tilde\phe\bigl(c_1(E-F),\ldots,c_r(E-F)\bigr)\;,
$$
where the $c_j(E-F)$ are the Chern classes of the virtual bundle~$E-F$, and 
$\tilde\phe\in\C[z_1,\ldots,z_r]$ is the unique polynomial such that 
$$
\phe=\tilde\phe(\sigma_1,\ldots,\sigma_r)\;,
$$
where $\sigma_1,\ldots,\sigma_r\in\C[z_1,\ldots,z_r]$ are the elementary symmetric functions
on $n$ variables. 

\newrem If $\Sigma(f)$ is finite, then the number of points in~$\Sigma(f)$, counted with respect
to a suitable multiplicity, is 
$$
{1\over\nu}[(\nu+1)^{n+1}-1]=\sum_{k=0}^{n}{n+1\choose k+1}\nu^k\;;
$$
see, e.g., [AT1].

\newrem Since the total Chern classes of $T\P^n$ in terms of the first Chern class of $N$
is given by $c(T\P^n)=\bigl(1-c_1(N)\bigr)^{n+1}$, the total Chern class of $T\P^n-N^{\otimes\nu}$ is given by
$$
c(T\P^n-N^{\otimes\nu})={\bigl(1-c_1(N)\bigr)^{n+1}\over 1+\nu c_1(N)}
=\bigl(1-c_1(N)\bigr)^{n+1}\sum_{j=0}^n (-1)^j\nu^j c_1(N)^j
$$
and thus the Chern classes of $T\P^n- N^{\otimes\nu}$ are given by
$$
c_j(T\P^n-N^{\otimes\nu})=(-1)^j\sum_{k=0}^j {n+1\choose k}\nu^{j-k}\,
c_1^j(N)\;.
$$
In particular,
$$
\int_{\P^n} c_n(T\P^n- N^{\otimes\nu})=\sum_{k=0}^n{n+1\choose k}\nu^{n-k}\quad\hbox{and}\quad \int_{\P^n} c_1^n(T\P^n- N^{\otimes\nu})=(\nu+1+n)^n\;.
\neweq\eqremd
$$
Arguing in a similar way we get that 
the Chern classes of $(T\P^n\oplus N)- N^{\otimes\nu}$ are given by
$$
c_j\bigl((T\P^n\oplus N)- N^{\otimes\nu}\bigr)=(-1)^j
\sum_{k=0}^j\left[{n+1\choose k}-{n+1\choose k-1}\right]\nu^{j-k}\,c_1^j(H)\;,
$$
with the convention ${n+1\choose -1}=0$. 
 
The usefulness of a theorem like Theorem~1 is of course related to the possibility of computing the residues involved, and indeed we have
explicit formulas for the (generic) case of point components, expressed in terms of
Grothendieck residues (the residues at components of positive dimension can be expressed by using Lehmann-Suwa theory; see [LS], [S, Chapter~IV] and [ABT1, (6.1)]). To state the formulas we recall three definitions.

\newdef Let $\ca O^n_{w^o}$ be the ring of germs of holomorphic functions at~$w^o\in\C^n$.
Take $g_1,\ldots, g_n\in\ca O^n_{w^o}$ such that $w^o$ is an isolated zero of $g=(g_1,\ldots,g_n)$. 
Then the {\sl Grothendieck residue} of $h\in\ca O^n_{w^o}$ along $g_1,\ldots,g_n$ is defined (see [H], [L, Section~5], [LS,
Section~4], [S,~pp.105--107]) by the formula
$$
\Res_{w^o}\left[ h(w)\, dw_1\wedge\cdots\wedge dw_n \atop g_1,\ldots,g_n\right]=
\left({1\over 2\pi i}\right)^n\int_{\Gamma}{h \over g_1 \cdots
g_n}\,dw_1\wedge\cdots\wedge dw_n\;,
$$
where $\Gamma=\{|g_j(w)|=\eps\mid j=1,\ldots,n\}$ for $\eps>0$ small enough, oriented so
that $d\arg(g_1)\wedge \cdots \wedge d\arg(g_n)>0$. 

\newdef Let $\phe\in\C[z_1,\ldots,z_n]$ be a homogeneous symmetric polynomial in $n$ variables, and $L\colon V\to V$ an endomorphism
of an $n$-dimensional complex vector space~$V$. Then we set
$$
\phe(L)=\phe(\lambda_1,\ldots,\lambda_n)\;,
$$
where $\lambda_1,\ldots,\lambda_n\in\C$ are the eigenvalues of~$L$. 

\newdef We shall say that a homogeneous polynomial self-map $F=(F_0,\ldots,F_n)\colon\C^{n+1}
\to\C^{n+1}$ of degree~$\nu+1$ {\sl induces} the meromorphic self-map $f\colon\P^n\dasharrow\P^n$
if 
$$
f([z_0:\cdots:z_n])=[F_0(z_0,\ldots,z_n):\cdots:F_n(z_0,\ldots,z_n)]
$$
for all $[z_0:\cdots:z_n]\in\P^n$. It is well known (see, e.g., [FS]) that every meromorphic self-map
of $\P^n$ of degree $\nu+1$ is induced by a unique homogeneous polynomial self-map
of~$\C^{n+1}$ of degree~$\nu+1$.

We can now state our next theorem:

\newthm Theorem \comput:  Let $f\colon\P^n\dasharrow\P^n$ be a meromorphic self-map of degree~$\nu+1\ge 2$ of the complex $n$-dimensional projective space. Let $\Sigma(f)=\Fix(f)\cup I(f)$ be
the union of the indeterminacy set $I(f)$ of~$f$ and the fixed points set $\Fix(f)$ of~$f$, and 
assume that $p=[1:w_1^o:\cdots:w_n^o]$ is an isolated point in~$\Sigma(f)$.  
Set $w^o=(w_1^o,\ldots,w_n^o)$, let $F=(F_0,\ldots,F_n)\colon\C^{n+1}\to\C^{n+1}$ be the homogeneous polynomial self-map of degree~$\nu+1$ inducing~$f$,
and define $g\colon\C^n\to\C^n$ by setting 
$$
g_j(w)=F_j(1,w)-w_j F_0(1,w)
\neweq\eqg
$$
for all $w\in\C^n$ and $j=1,\ldots,n$. Then:
\smallskip
\item{\rm(i)}we have
$$
\Res^1(f,p)= \Res_{w^o}\left[ F_0(1,w)^n\, dw_1\wedge\cdots\wedge dw_n \atop g_1,\ldots,g_n\right]\;;
$$
\item{\rm(ii)} if $\phe\in\C[z_1,\ldots,z_n]$ is a homogeneous symmetric polynomial of degree~$n$ then
$$
\Res^2_\phe(f,p)=\Res_{w^o}\left[ \phe(dg_{w})\, dw_1\wedge\cdots\wedge dw_n \atop g_1,\ldots,g_n\right]\;;
$$
\item{\rm(iii)} if $\psi\in\C[z_0,\ldots,z_n]$ is a homogeneous symmetric polynomial of degree~$n$ then
$$
\Res^3_\psi(f,p)=\Res_{w^o}\left[ \psi\bigl(F_0(1,w),\mu_1(w),\ldots,\mu_n(w)\bigr)\, dw_1\wedge\cdots\wedge dw_n \atop g_1,\ldots,g_n\right]\;,
$$
where $\mu_1(w),\ldots,\mu_n(w)$ are the eigenvalues of~$dg_w$.

This statement is effective because we can explicitly compute a Grothendieck residue using an algorithm suggested by Hartshorne (see [BB] and [H]). Without loss of generality, we can assume $w^o=O$.
Since the origin is an isolated zero of $g$, there exist minimal positive integers $\alpha_1,\ldots,
\alpha_n$ such that $w_1^{\alpha_1},\ldots,w_n^{\alpha_n}$ belong 
to the ideal generated by $g_1,\ldots,g_n$ in $\ca O^n=\ca O^n_O$. Hence there exist holomorphic functions $b_{ij}\in\ca O^n$ such that
$$
w_i^{\alpha_i}=\sum_{j=1}^n b_{ij} g_j\;.
$$
The properties of the Grothendieck residue (see [H]) then imply
$$
\Res_O\left[{h\,dw_1\wedge\cdots\wedge dw_n
\atop g_1,\cdots,g_n}\right]=\Res_O\left[
{h\det(b_{ij})\, dw_1\wedge\cdots\wedge dw_n\atop 
w_1^{\alpha_1},\cdots,w_n^{\alpha_n}}\right]\;.
$$
The right-hand side is now evaluated by expanding $h\det(b_{ij})$ in a
power series in the $w_i$; the residue is given by the coefficient of $w_1^{\alpha_1-1}\cdots w_n^{\alpha_n-1}$.

The easiest case is when $\det(dg_O)\ne 0$. Indeed, since $O$ is an isolated zero of $g$
we can write
$$
g_j=\sum_{i=1}^k c_{ji} w_i\;.
$$
Differentiating this and evaluating in~$O$ we get
$$
{\de g_j\over\de w_k}(O)=c_{jk}(O)\;;
$$
therefore if $\det(dg_{O})\ne 0$ we can invert the matrix $(c_{ji})$ in a neighborhood of~$O$ and
write
$$
w_i=\sum_{j=1}^n b_{ij} g_j\;,
$$
where $(b_{ij})$ is the inverse matrix of $(c_{ji})$. So we have $\alpha_1=\cdots=\alpha_n=1$
and $\det(b_{ij})(O)=1/\det(dg_{O})$. 

It thus follows that
$$
\det(dg_{w^o})\ne 0\quad\Longrightarrow\quad \Res_{w^o}\left[ h(w)\, dw_1\wedge\cdots\wedge dw_n \atop g_1,\ldots,g_n\right]=
{h(w^o)\over \mu_1\cdots\mu_n}\;,
\neweq\equzero
$$
where $\mu_1,\ldots,\mu_n\in\C$ are the eigenvalues of~$dg_{w^o}$.

In our situation, 
$g_j(w)=F_j(1,w)-w_j F_0(1,w)$;
therefore we have
$$
{\de g_j\over\de w_k}(w)={\de F_j\over\de w_k}(1,w)-w_j{\de F_0\over\de w_k}(1,w)-\delta^j_k F_0(1,w)\;,
$$
where $\delta^j_k$ is Kronecker's delta.

So if $p=[1:w_1^o:\cdots:w_n^o]$ is an indeterminacy point we have $F_0(1,w^o)=0$ and
$$
\Jac(g)(w^o)=\left({\de F_j\over\de w_k}(1,w^o)-w^o_j{\de F_0\over\de w_k}(1,w^o)\right)_{j,k=1,\ldots,n}\;.
$$
If instead $p$ is a fixed point, we can consider the 
differential of~$f$ at~$p$. Let $\chi$ be the usual chart of~$\P^n$ centered at $[1:0:\cdots:0]$. Then
$$
\tilde f(w)=\chi\circ f\circ\chi^{-1}(w)=\left({F_1(1,w)\over F_0(1,w)},\ldots,{F_n(1,w)\over
F_0(1,w)}\right)\;,
$$
and so
$$
{\de\tilde f_j\over\de w_k}(w)={1\over F_0(1,w)}\left[{\de F_j\over\de w_k}(1,w)-{F_j(1,w)\over
F_0(1,w)}{\de F_0\over\de w_k}(1,w)\right]\;.
$$
In particular,
$$
{\de\tilde f_h\over\de w_k}(w^o)={1\over F_0(1,w^o)}\left[{\de F_j\over\de w_k}(1,w^o)-w^o_j
{\de F_0\over\de w_k}(1,w^o)\right]\;;
\neweq\equdue
$$
therefore
$$
\Jac(g)(w^o)=F_0(1,w^o)\bigl[\Jac(\tilde f)(w^o)-I\bigr]\;.
$$
It follows that the eigenvalues $\mu_1,\ldots,\mu_n$ of $dg_{w^o}$ are related to the eigenvalues
$\lambda_1,\ldots,\lambda_n$ of~$df_p$ by the formula 
$$
\forevery{j=1,\ldots,n}\mu_j=F_0(1,w^o)(\lambda_j-1)\;.
$$
In particular, $\det(dg_{w^o})\ne 0$ if and only if 1 is not an eigenvalue of $df_p$, that is if and
only if $p$ is a simple fixed point of~$f$. 

Summing up, these computations give the following particular case of Theorem~\comput:

\newthm Corollary \computsimple:  Let $f\colon\P^n\dasharrow\P^n$ be a meromorphic self-map of degree~$\nu+1\ge 2$ of the complex $n$-dimensional projective space. 
\smallskip
\item{\rm(a)} Assume that $p\in\P^n$ is a simple (necessarily isolated) fixed point of $f$, and let 
$\lambda_1,\ldots,\lambda_n\ne 1$ be the eigenvalues of~$df_p$. Then:
\smallskip
\itemitem{\rm(i)}we have
$$
\Res^1(f,p)={(-1)^n\over(1-\lambda_1)\cdots(1-\lambda_n)}\;;
$$
\itemitem{\rm(ii)} if $\phe\in\C[z_1,\ldots,z_n]$ is a homogeneous symmetric polynomial of degree~$n$ then
$$
\Res^2_\phe(f,p)={\phe(1-\lambda_1,\ldots,1-\lambda_n)\over(1-\lambda_1)\cdots(1-\lambda_n)}\;;
$$
\itemitem{\rm(iii)} if $\psi\in\C[z_0,\ldots,z_n]$ is a homogeneous symmetric polynomial of degree~$n$ then
$$
\Res^3_\psi(f,p)=(-1)^n{\psi(1,\lambda_1-1,\ldots,\lambda_n-1)\over(1-\lambda_1)\cdots(1-\lambda_n)}\;.
$$
\item{\rm(b)} Assume that $p=[1:w^o_1:\cdots:w_n^o]$ is an isolated indeterminacy point, and 
that $\det G\ne 0$, where 
$$
G=\left({\de F_j\over\de w_k}(1,w^o)-w^o_j{\de F_0\over\de w_k}(1,w^o)\right)_{j,k=1,\ldots,n}\;,
$$
and $F=(F_0,\ldots,F_n)\colon\C^{n+1}\to\C^{n+1}$ is the homogeneous polynomial self-map of degree~$\nu+1$ inducing~$f$. Denote by $\mu_1,\ldots,\mu_n\ne 0$ the eigenvalues of~$G$. Then:
\itemitem{\rm(i)}we have
$$
\Res^1(f,p)=0\;;
$$
\itemitem{\rm(ii)} if $\phe\in\C[z_1,\ldots,z_n]$ is a homogeneous symmetric polynomial of degree~$n$ then
$$
\Res^2_\phe(f,p)={\phe(\mu_1,\ldots,\mu_n)\over\mu_1\cdots\mu_n}\;;
$$
\itemitem{\rm(iii)} if $\psi\in\C[z_0,\ldots,z_n]$ is a homogeneous symmetric polynomial of degree~$n$ then
$$
\Res^3_\psi(f,p)={\psi(0,\mu_1,\ldots,\mu_n)\over\mu_1\cdots\mu_n}\;.
$$

This corollary shows that 
our Theorem~\index\ is related to Ueda's index theorem, which however applies
only to holomorphic self-maps having only simple fixed points:

\newthm Theorem \Ueda: (Ueda [U]) Let $f\colon\P^n\to\P^n$ be holomorphic of degree~$\nu+1\ge 2$, and assume that all fixed points of $f$ are simple (and thus isolated). For $p\in\Fix(f)$, let 
$\lambda_1(p),\ldots,\lambda_n(p)\ne 1$ denote the eigenvalues of~$df_p$. Then
$$
\forevery{k=0,\ldots,n} \sum_{p\in\Fix(f)}{\sigma_k(df_p)\over\bigl(1-\lambda_1(p)\bigr)
\cdots\bigl(1-\lambda_n(p)\bigr)}=(-1)^k(\nu+1)^k\;.
$$

For instance, the case $k=0$ is a consequence of Theorem~\index.(i), and the cases
$k>0$ can be linked to Theorem~\index.(ii) and (iii) by using the formula
$$
\sigma_j(I-L)=\sum_{\ell=0}^j{n-\ell\choose n-j}(-1)^\ell \sigma_\ell(L)
$$
valid for every endomorphism $L$ of an $n$-dimensional complex vector space.
 
\newex
Let us consider
the map $f([z])=[z_0^{\nu+1}:\cdots:z_n^{\nu+1}]$. This map is holomorphic, and all its fixed
points are simple. More precisely, for each $\ell=0,\ldots,n$ the set $\Fix(f)$ contains exactly
${n+1\choose \ell+1}\nu^\ell$ fixed points of the form $[\zeta_0:\cdots:\zeta_n]$, where 
$\ell+1$ of the $\zeta_j$ are $\nu$-th roots of unity, and the remaining $\zeta_j$'s are equal
to~0. If $p\in\Fix(f)$ has exactly $\ell+1$ non-zero homogeneous coordinates, we shall say
that $p$ is a fixed point of {\sl level}~$\ell$.\hfill\break
\indent Using \equdue\ it is easy to see that if $p\in\Fix(f)$ is a fixed point of level~$\ell$ then
the eigenvalues of~$df_p$ are $\nu+1$ with multiplicity~$\ell$ and $0$ with multiplicity~$n-\ell$. 
In particular, $\bigl(1-\lambda_1(p)\bigr)\cdots\bigl(1-\lambda_n(p)\bigr)=(-1)^\ell\nu^\ell$,
and so Theorem~\index.(i) becomes
$$
(-1)^n\sum_{\ell=0}^n{n+1\choose\ell+1}(-1)^\ell=(-1)^n\;,
$$
while Theorem~\index.(ii) becomes
$$
(-1)^n\sum_{\ell=0}^n {n+1\choose\ell+1}(-1)^\ell\phe(\underbrace{\nu,\ldots,\nu}_{\ell\hbox{ \sevenrm times}},-1,\ldots,-1)
=\int_{\P^n}\phe(T\P^n- N^{\otimes\nu})\;,
$$
and Theorem~\index.(iii) becomes
$$
(-1)^n\sum_{\ell=0}^n {n+1\choose\ell+1}(-1)^\ell\psi(1,\underbrace{\nu,\ldots,\nu}_{\ell\hbox{ \sevenrm times}},-1,\ldots,-1)
=\int_{\P^n}\psi\bigl((T\P^n\oplus N)- N^{\otimes\nu}\bigr)\;.
$$
For instance, taking $\phe=\sigma_1^n$ and recalling \eqremd\
we get
$$
(\nu+1+n)^n=\int_{\P^n}c_1^n(T\P^n- N^{\otimes\nu})=(-1)^n\sum_{\ell=0}^n{n+1\choose\ell+1}(-1)^\ell
\bigl((\nu+1)\ell-n\bigr)^n\;.
$$
The (non-trivial) equality of the left and right-hand sides in this formula can also be proved directly
by using Abel's formula
$$
(x+y)^r=\sum_{k=0}^r {r\choose k}x(x-kz)^{k-1}(y+kz)^{r-k}
\neweq\eqAbel
$$ 
valid for all $r\in\N$ and $x$, $y$, $z\in\C$ with $x\ne 0$ (see, e.g., [K, p.~56]). Indeed we have
$$
\eqalign{
\sum_{\ell=0}^n{n+1\choose\ell+1}&(-1)^\ell\bigl((\nu+1)\ell-n\bigr)^n\cr
&=
\sum_{\ell=0}^n{n+1\choose\ell+1}\bigl(n-(\nu+1)\ell\bigr)^\ell\bigl(-n+(\nu+1)\ell\bigr)^{n-\ell}\cr
&=\sum_{k=1}^{n+1}{n+1\choose k}\bigl(n+\nu+1-k(\nu+1)\bigr)^{k-1}\bigl(-n-\nu-1+k(\nu+1)\bigr)^{n+1-k}\cr
&={1\over n+\nu+1}\biggl[\sum_{k=0}^{n+1}{n+1\choose k}(n+\nu+1)\bigl(n+\nu+1-k(\nu+1)\bigr)^{k-1}\bigl(-n-\nu-1+k(\nu+1)\bigr)^{n+1-k}\cr
&\phantom{={1\over n+\nu+1}}\quad+(-1)^n(n+\nu+1)^{n+1}\biggr]\cr
&=(-1)^n(\nu+n+1)^n\;,
\cr}
$$
thanks to \eqAbel\ applied with $r=n+1$, $x=-y=n+\nu+1$ and $z=\nu+1$.\hfil\break\indent
Analogously, taking $\phe=\sigma_n$ we get
$$
\sum_{k=0}^n{n+1\choose k}\nu^{n-k}=\int_{\P^n}c_n(T\P^n- N^{\otimes\nu})=\sum_{\ell=0}^{n}{n+1\choose\ell+1}\nu^\ell
={(\nu+1)^{n+1}-1\over\nu}=\sum_{\ell=0}^n(\nu+1)^\ell\;,
$$
in agreement with both \eqremd\ and Theorem~\Ueda.

Let us finally show how to use the local theory developed in [ABT1] and [AT2] to prove 
Theorems~\index\ and~\comput. 

Given a meromorphic $f\colon\P^n\dasharrow\P^n$ of degree $\nu+1\ge 2$, let $F=(F_0,\ldots,F_n)\colon
\C^{n+1}\to\C^{n+1}$ be the homogeneous polynomial map of degree~$\nu+1$ inducing~$f$.
First of all, we associate to $F$ the homogeneous vector field
$$
Q=\sum_{j=0}^n F_j{\de\over\de z_j}\;.
$$
It is well-known that the time-1 map $f_Q$ of~$Q$ is a germ of holomorphic self-map
of~$\C^{n+1}$ tangent to the identity; furthermore we can write $f_Q(z)=z+F(z)+O(\|z\|^{v+2})$.
In particular, by definition a direction $v\in\C^{n+1}\setminus\{O\}$ is 
a non-degenerate characteristic direction for $f_Q$ if and only if $[v]\in\P^n$ is a fixed
point of~$f$; and it is a degenerate characteristic direction for~$f_Q$ if and only if $[v]$ is
an indeterminacy point of~$f$. Furthermore, since $f$ has degree at least~2 then
$Q$ is non-dicritical, that is not all directions are characteristic.

Now let $\pi\colon M\to\C^{n+1}$ be the blow-up of the origin, and denote by $E=\pi^{-1}(O)$
the exceptional divisor. By construction, $E$ is canonically biholomorphic to~$\P^n$;
furthermore, the blow-up~$M$ can be identified with the total space of the normal bundle
$N_E$ of~$E$ in~$M$, and $N_E$ is isomorphic to the tautological line bundle $N$ on~$\P^n$.

We can now lift $f_Q$ to the blow-up, obtaining (see, e.g., [A]) a germ $\hat{f}_Q$ about~$E$ of 
holomorphic self-map of~$M$ fixing $E$ pointwise, and we may apply all the machinery
developed in~[ABT1]. First of all, since $Q$ is non-dicritical then $\hat{f}_Q$ is tangential
and has order of contact $\nu$ with~$E$ ([ABT1, Proposition~1.4]). 
Then we can define ([ABT1, Proposition~3.1 and Corollary~3.2]) the canonical section $X_f$, which is a global holomorphic section
of $TE\otimes(N_E^*)^{\otimes\nu}$. In local coordinates in $M$
centered at $[1:0:\cdots:0]\in E$ we can write
$$
X_f=\sum_{j=1}^n g_j\,{\de\over\de w^j}\otimes(dw_0)^{\otimes\nu}\;,
\neweq\eqXf
$$
where $g_1,\ldots,g_n$ are given by~\eqg, and $E$ in these coordinates is given by $\{w_0=0\}$.

In particular, we can think of $X_f$ as a bundle morphism $X_f\colon N_E^{\otimes\nu}\to TE$,
vanishing exactly on the characteristic directions of~$f_Q$, and thus we have proved
the following

\newthm Proposition \folia: Let $f\colon\P^n\dasharrow\P^n$ be a meromorphic self-map of degree
$\nu+1\ge 2$. Then the canonical section~$X_f$, given locally by $\eqXf$, defines
a global singular holomorphic foliation of~$\P^n$ in Riemann surfaces, singular exactly 
at the fixed and indeterminacy points of~$f$.

We can now apply the index theorems proved in [ABT1]. Theorem~\index.(i) follows
immediately from the Camacho-Sad-like index Theorem~6.2 in [ABT1], 
recalling that $\int_{\P^n}c_1(N)^n=(-1)^n$. Theorem~\index.(ii) follows from the
Baum-Bott index theorem (see [S, Th. I\negthinspace I\negthinspace I.7.6] and
[ABT1, Theorem~6.4]); and Theorem~\index.(iii) follows from the 
Lehmann-Suwa-like index Theorem~6.3 in [ABT1], that can be applied to this 
situation because $E$ is automatically comfortably embedded in~$M$ ([ABT1, Example~2.4]),
and recalling that the exact sequence
$$
O \longrightarrow TE \longrightarrow TM|_E \longrightarrow N_E \longrightarrow O
$$
implies that $c(TM|_E)=c(T\P^n\oplus N)$.

Theorem~\comput\ also follows from [ABT1]. Indeed, Theorem~\comput.(i) is an immediate
consequence of [ABT1, Theorem~6.5.(i)] and the computations in [AT2, Section~5].
The same computations and [ABT1, Theorem~6.6] yield Theorem~\comput.(iii); 
and Theorem~\comput.(ii) follows from [BB] and [S, Th. I\negthinspace I\negthinspace I.5.5].

\smallskip

{\it Acknowledgments.} I would like to thank Francesca Tovena, Jasmin Raissy and Matteo Ruggiero for many
useful conversations, Nur\`\i a Fagella and the Institut de Matem\`atica de la Universitat de Barcelona 
for their wonderful hospitality during the preparation of this note, and Jack Milnor for
creating so much beautiful mathematics and for being of inspiration to us all.

\setref{ABT3}
\beginsection References.

\art A M. Abate: Diagonalization of non-diagonalizable discrete
holomorphic dynamical systems! Amer. J. Math.! 122 2000 757-781

\art ABT1 M. Abate, F. Bracci, F. Tovena: Index theorems for holomorphic
self-maps! Ann. of Math.! 159 2004 819-864

\art ABT2 M. Abate, F. Bracci, F. Tovena: Index theorems for holomorphic
maps and foliations! Indiana Univ. Math. J.! 57 2008 2999-3048

\art ABT3 M. Abate, F. Bracci, F. Tovena: Embeddings of submanifolds and
normal bundles! Adv. Math.! 220 2009 620-656

\art AT1 M. Abate, F. Tovena: Parabolic curves in $\C^3$! Abstr. Appl. 
Anal.! 2003 2003 275-294

\pre AT2 M. Abate, F. Tovena: Poincar\'e-Bendixson theorems for meromorphic connections and homogeneous vector fields! To appear in J. Diff. Equations. ArXiv:0903.3485! 2011

\art BB P. Baum and R. Bott: Singularities of holomorphic
foliations! J. Diff. Geom.! 7 1972 279-342


\art FS J.-E. Forn\ae ss, N. Sibony: Complex dynamics in higher dimension. I! Ast\'erisque! 222 1994 201-231

\book H R. Hartshorne: Residues and duality! Lecture Notes in Math. Vol. 20, Springer,
Berlin, 1966

\book K D.E. Knuth: Fundamental algorithms! Addison-Wesley, Reading, MA, 1973

\art L D. Lehmann: R\'esidues des sous-vari\'et\'es invariants
d'un feuilletage singulier! Ann. Inst. Fourier, Grenoble! 41 1991
211-258

\art LS D. Lehmann and T. Suwa: Residues of holomorphic vector
fields relative to singular invariant subvarieties! J. Diff.
Geom.! 42 1995 165-192

\book Su T. Suwa: Indices of vector fields and residues of
singular holomorphic foliations! Hermann, Paris, 1998

\coll U T. Ueda: Complex dynamics on projective spaces---index formula for fixed points!
Dynamical systems and chaos, Vol. 1 (Hachioji, 1994)! World Sci. Publ., River Edge, NJ, 1995,
pp. 252--259

\bye

*****



Ueda's index theorem is a result valid for holomorphic self-maps of a projective space; but we shall see that it is related to results in local holomorphic dynamics too.

Let us introduce some notations. We shall denote by $\P^n$ the complex projective space 
of dimension~$n$, and by $\pi\colon\C^{n+1}\setminus\{O\}\to\P^n$ the standard projection
$$
(z^0,\ldots,z^n)\mapsto \pi(z^0,\ldots,z^n)=[z^0:\cdots:z^n]\;,
$$
where $[z^0:\cdots:z^n]$ are homogeneous coordinates. We shall sometimes write $[z]$ instead of~$\pi(z)$.

Furthermore, we shall denote by $\widehat{\C^{n+1}}$ the blow-up of the origin in$\C^{n+1}$,
and by $\hat\pi\colon\widehat{\C^{n+1}}\to\C^{n+1}$ the canonical projection. If $E=\hat\pi^{-1}(O)$
is the exceptional divisor, it is well-known that $E$ is the projective space of complex tangent directions at the origin, and thus it can be identified with $\P^n$. Furthermore, $\widehat{\C^{n+1}}$
can be canonically identified with the total space of the normal bundle $\pi\colon N_E\to E$
of~$E$ in $\widehat{\C^{n+1}}$ in such a way that $\C^{n+1}\setminus\{O\}\equiv N_E\setminus E$
(where we are identifying $E$ with the zero section of~$N_E$). Finally, under these identifications
the canonical projection of $\C^{n+1}\setminus\{O\}$ onto~$\P^n$ coincides with the
projection of $N_E\setminus E$ onto~$E$, thus justifying the usage of the same symbol for both
projections.

There is a well-known relationships between holomorphic (or meromorphic) self-maps of
$\P^n$ and homogeneous self-maps of $\C^{n+1}$. To each meromorphic self-map 
$f\colon\P^n\to\P^n$ is associated a unique homogeneous polynomial map $F\colon\C^{n+1}\to\C^{n+1}$ such that $f\circ\pi=\pi\circ F$; furthermore, $f$ is holomorphic if and only if $F^{-1}(O)=\{O\}$. Using homogeneous coordinates, we are saying that every meromorphic $f\colon\P^n\to\P^n$
is of the form
$$
f([z^0:\cdots:z^n])=[F^0(z):\cdots:F^n(z)]\;,
\neweq\eqduno
$$
where $z=(z^0,\ldots,z^n)$ and $F^0,\ldots,F^n$ are homogeneous polynomials of the same degree~$\nu$, called the {\sl degree} of~$f$. Conversely, given a homogeneous polynomial map
$F=(F^0,\ldots,F^n)\colon\C^{n+1}\to\C^{n+1}$ of degree~$\nu$, formula \eqduno\ defines 
a meromorphic self-map of~$\P^n$, holomorphic if and only if $F^0(z)=\cdots=F^n(z)=0$ only if $z=O$.

This double interpretation implies that in the literature different terminologies have been introduced. A {\sl characteristic direction} for a homogeneous polynomial map $F\colon\C^{n+1}\to\C^{n+1}$ is a non-zero vector $v\in\C^{n+1}\setminus\{O\}$
such that $F(v)=\lambda v$ for a suitable $\lambda\in\C$. The line $L_v=\C v$ will be a
{\sl characteristic line} for $F$. Moreover, we shall say that
$v$ is {\sl non-degenerate} if $\lambda\ne 0$, and {\sl degenerate} otherwise.

If $f\colon\P^n\to\P^n$ is the meromorphic self-map associated to $F\colon\C^{n+1}\to\C^{n+1}$,
it is clear that $v\in\C^{n+1}\setminus\{O\}$ is a non-degenerate characteristic direction
for $F$ if and only if $[v]\in\P^n$ is a fixed point for $f$; and that it is a degenerate characteristic direction
for~$F$ if and only if $[v]$ is an indeterminacy point for $f$.

Furthermore, to each fixed point $[v]$ for $f$ correspond $\nu-1$ distinct fixed points for $F$ in
the characteristic line $L_v$: indeed, if $v\in\C^{n+1}\setminus\{O\}$ is any representative
of $[v]$, so that $F(v)=\lambda v$, we have $F(\alpha v)=\alpha^\nu\lambda v$, which is
equal to $\alpha v$ if and only if $\alpha$ is a $(\nu-1)$-root of $\lambda^{-1}$. 

We can also define a multiplicity for isolated characteristic directions (or fixed or indeterminacy points).  If
$v$ is an isolated characteristic direction of $F$, its {\sl multiplicity}~$\mu_F(v)$ is the
local intersection multiplicity (see, e.g., [C] or [GH] for definition and properties of
the local intersection multiplicity) at~$v$ in~$\P^{n-1}$ of the polynomials
$z_{j_0}P_j-z_jP_{j_0}$ with $j\ne j_0$, where $j_0$ is any index such that $v_{j_0}\ne 0$
(and $\mu_P(v)$ is clearly independent of $j_0$).

\newthm Lemma \FS: Let $F=(F^0,\ldots,F^n)\in\End(\C^{n+1},O)$ be a homogeneous
map of degree $\nu\ge 2$. Denote by
$\pi\colon\P^{n+1}\setminus\{[0:\cdots:0:1]\}\to\P^n$ the projection
$\pi([v^0:v^1:\cdots:v^{n+1}])=[v^0:\cdots:v^n]$, and let~$\ca S\subset\P^{n+1}$ be the set of
solutions of the system 
$$
\cases{P^0(z)-(z^{n+1})^{\nu-1}z^0=0,\cr 
\hfill\vdots\hfill\cr 
P^n(z)-(z^{n+1})^{\nu-1}z^n=0.}
\neweq\eqFS
$$
Then
\smallskip
\item{\rm(i)} The vector $v=(v^0,\ldots,v^n)\in\C^{n+1}\setminus\{O\}$ is a characteristic direction for
$F$ iff $\pi^{-1}([v])\cap\ca S$ is not empty.
More precisely, $v$ is a degenerate characteristic direction iff $\pi^{-1}([v])\cap\ca
S=\{[v^0:\cdots:v^n:0]\}$, and it is a non-degenerate characteristic
direction iff $\pi^{-1}([v])\cap\ca S$ contains exactly $\nu-1$ elements all with non-zero
last coordinate.
\item{\rm(ii)} If $v$ is a non-degenerate isolated characteristic direction, then its
multiplicity~$\mu_F(v)$ is equal to the multiplicity  of any element in
$\pi^{-1}([v])\cap\ca S$ as solution of $\eqFS$; on the other hand, if $v$ is a degenerate
isolated characteristic direction, then the multiplicity of $[v^0:\cdots:v^n:0]$ as
solution of~$\eqFS$ is~$(\nu-1)\mu_F(v)$. 
\item{\rm(iii)} The number of characteristic directions of $F$, counted according to
their multiplicities, if finite is given\break\indent by $(\nu^n-1)/(\nu-1)$.

\pf (i) This is obvious.

(ii) Without loss of generality we can assume that $v=[0:\cdots:0:1]$, and fix $\tilde
v\in\pi^{-1}(v)\cap\ca S$. In the usual local coordinates of the subset $\{z_n\ne
0\}\subset\P^n$ the point~$\tilde v$ is represented by $(\lambda,0,\ldots,0)$, where
$\lambda=0$ iff $v$ is degenerate. Analogously, the local coordinates of the subset $\{z_n\ne
0\}\subset\P^{n-1}$ are centered in~$v$. So we have
$$
\mu_P(v)=I\bigl(P_1(z',1)-z_1P_n(z',1),\ldots,P_{n-1}(z',1)-z_{n-1}P_n(z',1); O\bigr),
$$
while the multiplicity of $\tilde v$ as solution of \eqFS\ is
$$
\tilde\mu=I\bigl(P_1(z',1)-z_0^{\nu-1}z_1,\ldots,P_{n-1}(z',1)-z_0^{\nu-1}z_{n-1},
	P_n(z',1)-z_0^{\nu-1};(\lambda,0,\ldots,0)\bigr),
$$
where $I$ denotes the local intersection multiplicity, and $z'=(z_1,\ldots,z_{n-1})$. The
standard properties of
$I$ immediately yields
$$
\tilde\mu=I\bigl(P_1(z',1)-z_1P_n(z',1),\ldots,P_{n-1}(z',1)-z_{n-1}P_n(z',1),
P_n(z',1)-z_0^{\nu-1}; (\lambda,0,\ldots,0)\bigr).
$$
Set $Q_j(z')=P_j(z',1)-z_jP_n(z',1)$ for $j=1,\ldots,n-1$. Now, if $v$ is degenerate, that is
if $\lambda=0$, the ring $\ca O_n/\bigl(Q_1,\ldots,Q_{n-1},P_n(z',1)-z_0^{\nu-1}\bigr)$
is generated by $1,z_0,\ldots,z_0^{\nu-2}$ on the ring $\ca
O_{n-1}/\bigl(Q_1,\ldots,Q_{n-1}\bigr)$; therefore
$$
\tilde\mu=(\nu-1)\mu_P(v),
$$
as claimed.

On the other hand, if $\lambda\ne 0$ we can translate the coordinates obtaining
$$
\tilde\mu=I\bigl(Q_1,\ldots,Q_{n-1},P_n(z',1)-(z_0+\lambda)^{\nu-1};O\bigr).
$$
Now, $P_n(O',1)=\lambda^{\nu-1}$; therefore there is $Q_n\in\ca O_{n-1}$ such that
$Q_n(O')=\lambda$ and
$$
P_n(z',1)-(z_0+\lambda)^{\nu-1}=\bigl(Q_n(z')-z_0-\lambda\bigr)R(z_0,z')
$$
in $\ca O_n$, with $R(O)\ne 0$. Indeed, it suffices to take $Q_n\in\ca O_{n-1}$ so that
$Q_n(O')=\lambda$ and $Q_n^{\nu-1}=P_n(z',1)$ in~$\ca O_{n-1}$. But then
$$
\tilde\mu=I\bigl( Q_1,\ldots,Q_{n-1},Q_n-\lambda-z_0; O\bigr),
$$
and thus the argument used before yields $\tilde\mu=\mu_P(v)$.

(iii) By Bezout's theorem we know that \eqFS\ has exactly (infinite or) $\nu^n$
solutions, counted according to their multiplicities. The solution $[1:0:\cdots:0]$, which is
the only one not generating a characteristic direction of~$P$, has multiplicity 1. Thus we
are left with
$\nu^n-1$ solutions, and the assertion follows from (i) and (ii).\qedn

There is another, equivalent, definition of multiplicity. Assume that $O\in\C^n$ is an isolated 
fixed point of a germ $f=(f^1,\ldots,f^n)\in\ca O^n$. Then the multiplicity of $O$ is the
local intersection multiplicity of $f^1-z^1,\ldots,f^n-z^n$ at~$O$. *******

A few useful formulas:

Assume $F(a,O')=(a,O')$, with $a\ne 0$. This means that
$$
F^0(1,O')=a^{1-\nu}\quad\hbox{and}\quad F^h(1,O')=0\quad\hbox{for $h=1,\ldots,n$.}
\neweq\equuno
$$
Furthermore,
$$
{\de F^h\over\de z^k}(a,O')=a^{\nu-1}{\de F^h\over\de z^k}(1,O')\;.
$$
Euler's formula 
$$
\sum_{j=0}^n z^j{\de F^h\over\de z^j}=\nu F^h
$$
yields
$$
{\de F^h\over\de z^0}(a,O')=\cases{\nu & if $h=0$,\cr 0 & if $h=1,\ldots,n$.\cr}
$$
Now let $\phe_0$ be the usual chart of~$\P^n$ centered at $[1:0:\cdots:0]=\pi(a,O')$. Then
$$
\tilde f(w)=\phe_0\circ f\circ\phe_0^{-1}(w)=\left({F^1(1,w)\over F^0(1,w)},\ldots,{F^n(1,w)\over
F^0(1,w)}\right)\;,
$$
and so
$$
{\de\tilde f^h\over\de w^k}(w)={1\over F^0(1,w)}\left[{\de F^h\over\de z^k}(1,w)-{F^h(1,w)\over
F^0(1,w)}{\de F^0\over\de z^k}(1,w)\right]\;.
$$
In particular,
$$
{\de\tilde f^h\over\de w^k}(O)=a^{\nu-1}{\de F^h\over\de z^k}(1,O')={\de F^h\over\de z^k}(a,O')\;.
\neweq\equdue
$$
Summing up we get
$$
\hbox{Jac}(F)(a,O')=\pmatrix{\nu&\vrule& {\de F^0\over\de z^h}(a,O')\cr
\noalign{\hrule}
O&\vrule&\hbox{Jac}(\tilde f)(O)\cr}\;.
$$
Now, given a linear operator $L\colon V\to V$ of an $n$-dimensional $\C$-vector space~$V$,
we define $\sigma_1(L),\ldots,\sigma_n(L)\in\C$ by setting
$$
\det(I+tL)=1+t\sigma_1(L)+\cdots+t^n\sigma_n(L)=(-t)^n\det\bigl((-t)^{-1}I-L\bigr)
$$
or, equivalently (putting $s=(-t)^{-1}$)
$$
\det(sI-L)=s^n-s^{n-1}\sigma_1(L)+\cdots+(-1)^n\sigma_n(L)\;.
$$
If $\lambda_1,\ldots,\lambda_n\in\C$ are the eigenvalues of $L$ then $\det(sI-L)=(s-\lambda_1)
\cdots(s-\lambda_n)$; therefore $\sigma_j(L)$ is the $j$-th elementary simmetric function 
computed on the eigenvalues of~$L$.

If we set $\Phi_F(s,v)=\det(sI-dF_v)$ and $\Phi_f(s,[v])=\det(sI-df_{[v]})$, where $v$ is a fixed point
of~$F$ (and thus $[v]$ is a fixed point of~$f$), then
$$
\Phi_F(s,v)=(s-\nu)\Phi_f(s,[v])\;.
$$

If $[v]\in\P^n$ is a fixed point of $f$, a {\sl director} for $f$ at~$[v]$ is an eigenvalue of the
linear operator ${1\over\nu-1}(df_{[v]}-I)$. In other words, setting
$$
\mu={1\over\nu-1}(\lambda-1)\quad\Leftrightarrow\quad \lambda=(\nu-1)\mu+1
$$
we have that $\lambda$ is an eigenvalue of~$df_{[v]}$ if and only if $\mu$ is a director.
In particular, $[v]$ is simple (that is multiplicity~1) if and only if no director is zero. In particular,
$$
\Phi_f(1,[v])=(1-\nu)^n\mu_1\cdots\mu_n\;.
$$

Ueda's index theorem is a consequence of the

\newthm Theorem \Uedauno: Let $F\in\Hol(\C^{n+1},\C^{n+1}$ be homogeneous of degree~$\nu\ge 2$. Assume that all fixed points of $F$ are simple. Then
$$
\sum_{p\in\Fix(F)} {\Phi_F(t,p)\over\Phi_F(1,p)}=\nu^{n+1}\;.
\neweq\eqUedauno
$$

Then

\newthm Theorem \Uedadue: Assume that all fixed points of $f\in\Hol(\P^n,\P^n)$ are simple. Then
$$
\sum_{p\in\Fix(f)} {\Phi_f(t,p)\over\Phi_f(1,p)}={t^{n+1}-\nu^{n+1}\over t-\nu}=\sum_{k=0}^n
\nu^{n-k}t^k\;,
\neweq\eqUedadue
$$
where $\nu\ge 2$ is the degree of $f$.

\pf Since all the derivatives of $F$ vanish at the origin, $O$ is a simple fixed point of $F$ with
$\Phi_F(t,O)=t^{n+1}$. Therefore
$$
\sum_{p\in\Fix(F)} {\Phi_F(t,p)\over\Phi_F(1,p)}=\Phi_F(t,O)+(\nu-1)\sum_{p\in\Fix(f)}
{(t-\nu)\Phi_f(t,p)\over(1-\nu)\Phi_f(1,p)}=t^{n+1}-(t-\nu)\sum_{p\in\Fix(f)} {\Phi_f(t,p)\over\Phi_f(1,p)}\;,
$$ 
and the assertion follows from \eqUedauno.\qedn

\newthm Corollary \Uedatre: Assume that all fixed points of $f\in\Hol(\P^n,\P^n)$ are simple,
and if $p\in\Fix(f)$ then denote by $\lambda_1(p),\ldots,\lambda_n(p)$ the eigenvalues
of $df_p$. Then
$$
\sum_{p\in\Fix(f)}\left({1\over 1-\lambda_1(p)}+\cdots+{1\over 1-\lambda_n(p)}\right)=
\sum_{k=0}^n k\nu^{n-k}\;.
$$

\pf It suffices to differentiate both sides of \eqUedadue\ with respect to~$t$ and then let $t=1$.\qedn

\newthm Corollary \Uedaqua: Let $f\in\Hol(\P^n,\P^n)$. Then there is at least one fixed point $p$ of $f$
such that $df_p$ has an eigenvalue equal to~1 or with absolute value greater than~1.

\pf Suppose not. Then all fixed points of $f$ are simple and $|\lambda_j(p)|\le 1$ for all 
$j=1,\ldots,n$ and $p\in\Fix(f)$. Now, it is easy to check that
$$
\Re{1\over 1-\lambda}\ge {1\over 2}\quad\Leftrightarrow\quad |\lambda|^2\le 1\;;
$$
therefore 
$$
\sum_{p\in\Fix(f)}\left({1\over 1-\lambda_1(p)}+\cdots+{1\over 1-\lambda_n(p)}\right)
\ge {n\over 2}\sum_{k=0}^n \nu^k\;.
$$
But
$$
{n\over 2}\sum_{k=0}^n \nu^k-\sum_{k=0}^n k\nu^{n-k}=\sum_{k=0}^n\left({n\over2}-k\right)\nu^{n-k}
={n\over2}(\nu^n-1)+\left({n\over2}-1\right)(\nu^{n-1}-\nu)+\cdots>0\;,
$$
against Corollary~\Uedatre.\qedn

\newthm Corollary \Uedacin: Assume that all fixed points of $f\in\Hol(\P^n,\P^n)$ are simple. 
Then there are infinitely many periodic points of $f$ having an eigenvalue equal to~1 or with absolute value greater than~1.

\pf Let $S\subset\N$ be the set of natural numbers $r\ge 2$ such that $\lambda_j(p)^r\ne 1$
for all $j=1,\ldots,n $ and all $p\in\Fix(f)$. Clearly, $S$ contains arbitrarily large prime numbers;
moreover, all $p\in\Fix(f)$ are simple fixed points for $f^r$ if and only if $r\in S$.

We claim that when $r\in S$ is large enough, there exists $p\in\Fix(f^r)\setminus\Fix(f)$ 
such that $d(f^r)_p$ has an eigenvalue equal to 1 or with absolute value greater than 1.
Applying this claim to arbitrarily large prime numbers we get the assertion.

Assume then that all fixed points of $f^r$ are simple; then Corollary~\Uedatre\ yields
$$
\sum_{p\in\Fix(f)}\sum_{j=1}^n{1\over 1-\lambda_j(p)^r}+\sum_{q\in\Fix(f^r)\setminus\Fix(f)}
\sum_{j=1}^n{1\over 1-\mu_j(q)}=\sum_{k=0}^n k\nu^{r(n-k)}\;,
\neweq\eqUedacin
$$
where $\mu_1(q),\ldots,\mu_n(q)$ are the eigenvalues of $d(f^r)_q$ for $q\in\Fix(f^r)\setminus\Fix(f)$. 

Now, $\Re\bigl(1/(1-\lambda)\bigr)=1/2$ if and only if $|\lambda|=1$; and
$$
\lim_{r\to+\infty}\Re{1\over 1-\lambda^r}=
\cases{1 & if $|\lambda|<1$,\cr 0 & if $|\lambda|>1$.\cr}
$$
Therefore we can find a constant $K>0$ such that
$$
\Re\sum_{p\in\Fix(f)}\sum_{j=1}^n{1\over 1-\lambda_j(p)^r}\ge -K
$$
for all $r\in S$. 

Now assume, by contradiction, that $|\mu_j(q)|\le 1$ for all $j=1,\ldots,n$ and $q\in\Fix(f^r)\setminus
\Fix(f)$. Then
$$
\Re\sum_{q\in\Fix(f^r)\setminus\Fix(f)}
\sum_{j=1}^n{1\over 1-\mu_j(q)}\ge {n\over 2}\left({\nu^{r(n+1)}-1\over\nu^r-1}-
{\nu^{n+1}-1\over\nu-1}\right)\;.
$$
The right-hand side is a polynomial in~$\nu$ with leading term $(n/2)\nu^{nr}$. On the other
hand, the right-hand side of \eqUedacin\ is a polynomial in~$\nu$ with leading term $\nu^{r(n-1)}$;
therefore \eqUedacin\ cannot hold for $r$ large enough, and we are done.\qedn

\smallsect 3. Index theorems and partial holomorphic connections

There is another approach. In [AT2] we have shown that given $F\in\Hol(\C^{n+1},\C^{n+1})$
non-dicritical we can define a morphism $X_F\colon N_E^{\otimes(\nu-1)}\to TE$,
whose singular points are exactly the characteristic directions of $F$. In this way we get 
a rank-one foliation of $E=\P^n$, singular at the characteristic directions of $F$ (that is,
at the fixed or indeterminacy points of $f$). In the usual chart centered at $[1:0:\cdots:0]$
the morphism $X_F$ is given by
$$
X_F=\sum_{j=1}^n \bigl(F^j(1,w)-w^j F^0(1,w)\bigr){\de\over\de w^j}\otimes(d z^0)^{\otimes(\nu-1)}\;.
$$
We can then apply directly Baum-Bott theorem:

\newthm Theorem \BB: Let $M$ be a compact $n$-dimensional complex manifold, $L$ a line bundle over $M$, 
and $X\colon L\to TM$ be a holomorphic bundle map with isolated zeroes. 
Given a symmetric homogeneous polynomial $\phe$ of degree $n$, for each $p\in\Sing(X)$ 
we can define $\Res_\phe(X,p)\in\C$ depending only on the local behavior of $X$ at~$p$ so that 
$$
\sum_{p\in\Sing(X)}\Res_\phe(X,p)=\int_M \phe(TM-L)\;.
$$

Here, if $\phe\in\C[z_1,\ldots,z_n]$ is a symmetric homogeneous polynomial of degree~$n$, there is a unique
polynomial~$\tilde\phe\in\C[z_1,\ldots,z_n]$ such that 
$$
\phe=\tilde\phe(\sigma_1,\ldots,\sigma_n)\;,
$$
where $\sigma_1,\ldots,\sigma_n\in\C[z_1,\ldots,z_n]$ are the elementary symmetric functions
on $n$ indeterminates. So if $c_1(TM-L),\ldots,c_n(TM-L)$ are the Chern classes of the virtual
bundle $TM-L$ by definition we put
$$
\phe(TM-L)=\tilde\phe\bigl(c_1(TM-L),\ldots,c_n(TM-L)\bigr)\;.
$$
Analogously, if $L\colon V\to V$ is a complex linear operator of an $n$-dimensional vector space,
we put
$$
\phe(L)=\tilde\phe\bigl(\sigma_1(L),\ldots,\sigma_n(L)\bigr)\;.
$$
If $p\in\Sing(X)$, choose a local non-vanishing section $s$ of~$L$ near~$p$, a
system of coordinates $(w^1,\ldots,w^n)$ centered at~$p$, and write
$$
X(s)=\sum_{j=0}^n g^j {\de\over\de w^j}\;.
$$
Let $\Jac(g)=\bigl(\de g^h/\de w^k\bigr)$; then $\Res_\phe(X,p)$ is defined by the
Grothendieck residue
$$
\Res_\phe(X,p)=\Res_p\left[{\phe\bigl(\Jac(g)\bigr)\,dw^1\wedge\cdots\wedge dw^n
\atop g^1,\cdots,g^n}\right]\;.
$$
Let $\mu_1,\ldots,\mu_n$ be the eigenvalues of $\Jac(g)(p)$. If $\mu_1\cdots\mu_n\ne 0$, 
then the general properties of the Grothendieck residue given in [H] imply that
$$
\Res_\phe(X,p)={\phe(\mu_1,\ldots,\mu_n)\over\mu_1\cdots\mu_n}\;.
$$
More generally, Hartshorne (see [BB]) has given the following algorithm for computing
a Grothendieck residue. 
Since the origin is an isolated zero of the $g^j$, there exist minimal positive integers $\alpha_1,\ldots,
\alpha_n$ such that $(w^1)^{\alpha_1},\ldots,(w^n)^{\alpha_n}$ belong 
to the ideal generated by $g^1,\ldots,g^n$. Hence there exist holomorphic functions $b^i_j$	near $p$ such that
$$
(w^i)^{\alpha_i}=\sum_{j=1}^n b^i_j g^j\;.
$$
One then has
$$
\Res_p\left[{\phe\bigl(\Jac(g)\bigr)\,dw^1\wedge\cdots\wedge dw^n
\atop g^1,\cdots,g^n}\right]=\Res_p\left[
{\phe\bigl(\Jac(g)\bigr)\det(b^i_j)\, dw^1\wedge\cdots\wedge dw^n\atop 
(w^1)^{\alpha_1},\cdots,(w^n)^{\alpha_n}}\right]\;.
$$
The right hand side is now evaluated by expanding $\phe\bigl(\Jac(g)\bigr)\det(b^i_j)$ in a
power series in the $w^i$. The desired answer is then the coefficient of $(w^1\cdots w^n)^{-1}$ in	the resulting
Laurent series for $\phe\bigl(\Jac(g)\bigr)\det(b^i_j)/\bigl((w^1)^{\alpha_1}\cdots(w^n)^{\alpha_n}\bigr)$.

\newrem The easiest case is when $\det\Jac(g)(p)\ne 0$. Indeed, since $O$ is a singular point
of $X_F$, we can write
$$
g^j=\sum_{i=1}^k c^j_i w^i\;.
$$
Differentiating this and computing in~$p$ we get
$$
{\de g^j\over\de w^k}(p)=c^j_k(p)\;;
$$
therefore if $\det\Jac(g)(p)\ne 0$ we can invert the matrix $(c^j_i)$ in a neighborhood of~$p$ and
write
$$
w^i=\sum_{j=1}^n b^i_j g^j\;,
$$
where $(b^i_j)$ is the inverse matrix of $(c^i_j)$. So we have $\alpha_1=\cdots=\alpha_n=1$
and $\det(b^i_j)(p)=\det(c^i_j)(p)^{-1}=(\mu_1\cdots\mu_n)^{-1}$; thus it follows that
$$
\Res_\phe(X,p)={\phe(\mu_1,\ldots,\mu_n)\over\mu_1\cdots\mu_n}\;.
$$

In our situation, 
$$
g^j(w)=F^j(1,w)-w^j F^0(1,w)\;;
$$
therefore
$$
{\de g^j\over\de w^k}={\de F^j\over\de w^k}(1,w)-w^j{\de F^0\over\de w^k}(1,w)-\delta^j_k F^0(1,w)
$$
and
$$
{\de g^j\over\de w^k}(p)={\de F^j\over\de w^k}(1,O')-\delta^j_k F^0(1,O')\;.
$$
Now we have two cases to consider. 
\smallskip
\item{(a)} if $p$ is a degenerate characteristic direction, that is $[p]$ is an indeterminacy point of~$f$, then $F^0(1,O')=0$; 
\item{(b)} if $p$ is a non-degenerate characteristic direction, that is $[p]$ is a fixed point of $f$,
then $F^0(1,O')\ne 0$ and \equuno\ and \equdue\ yields
$$
\Jac(g)(p)=F^0(1,O') \bigl[\Jac(\tilde f)([p])-I\bigr]\;.
$$
In particular, the relation between the eigenvalues $\mu_j$ of $\Jac(g)$ and the eigenvalues
$\lambda_j$ of $df_p$ is
$$
\mu_j= F^0(1,O')(\lambda_j-1)\;.
$$
In particular, having $\det\Jac(g)\ne 0$ is equivalent to saying that $[p]$ is a simple fixed point of $f$, and in this case
$$
\Res_\phe(X_F,p)={\phe(\lambda_1-1,\ldots,\lambda_n-1)\over(\lambda_1-1)\cdots(\lambda_n-1)}\;.
$$

For the Camacho-Sad case, we recall that $h^1(w)=F^0(1,w)$ and that

\smallsect 4. Saito index formula

\smallsect 5. \'Ecalle's approach

\bye